\documentclass[12pt]{article}
\usepackage{color}
\usepackage{euscript,amssymb,amsfonts,amsmath,amsthm}
\usepackage{mathrsfs}
\usepackage{mathtext}
\usepackage{mathrsfs}
\usepackage[cp1251]{inputenc}
\usepackage[T2A]{fontenc}
\usepackage[english,russian]{babel}

\renewcommand{\P}{{\sf P}}
\newcommand{\E}{{\sf E}\,}

\newcommand{\eps}{\varepsilon}

\renewcommand{\kappa}{\varkappa}
\renewcommand{\le}{\leqslant}
\renewcommand{\ge}{\geqslant}

\newcommand{\R}{\mathbb R}

\newcommand{\F}{\mathcal F}
\newcommand{\ud}{\Delta_n} 

\newcommand{\sign}{\mathrm{sign}}

\newtheorem{theorem}{\sc Theorem}
\newtheorem{lemma}{\sc Lemma}
\newtheorem{corollary}{\sc Corollary}
\theoremstyle{remark}
\newtheorem{remark}{\sc Remark}

\renewcommand{\proof}{{\sc Proof}}

\title{On the absolute constants in the Berry--Esseen type
inequalities for identically distributed summands}
\author{Irina Shevtsova\footnote{Faculty
of Computational Mathematics and Cybernetics, Moscow State
University; Institute for Informatics Problems, Russian Academy of
Sciences; ishevtsova@cs.msu.su}}

\date{}

\begin{document}

\maketitle

{\small

\centerline{\bf Abstract}

\smallskip

By a modification of the method that was applied in (Korolev and
Shevtsova, 2010), here the inequalities
$$
\Delta_n\le 0.3328\cdot\frac{\beta_3+0.429}{\sqrt{n}},
$$
$$
\Delta_n\le 0.33554\cdot\frac{\beta_3+0.415}{\sqrt{n}}
$$
are proved for the uniform distance $\Delta_n$ between the standard
normal distribution function and the distribution function of the
normalized sum of an arbitrary number $n\ge1$ of independent
identically distributed random variables with zero mean, unit
variance and finite third absolute moment $\beta_3$. The first of
these two inequalities improves one that was proved in (Korolev and
Shevtsova, 2010), and as well sharpens the best known upper estimate
for the absolute constant $C_0$ in the classical Berry--Esseen
inequality to be $C_0<0.4756$, since
$0.3328(\beta_3+0.429)\le0.3328\cdot1.429\beta_3<0.4756\beta_3$ by
virtue of the condition $\beta_3\ge1$. The second of these
inequalities is also a structural improvement of the classical
Berry--Esseen inequality, and as well sharpens the upper estimate
for $C_0$ still more to be $C_0<0.4748$.

\smallskip


{\bf MSC classes:}  60F05

}

\section{Introduction and formulation of the main results}

By $\F_3$ we will denote the set of distribution functions with zero
mean, unit variance and finite third absolute moment $\beta_3$. Let
$X_1,X_2,\ldots$~be independent random variables with common
distribution function $F\in\F_3$ defined on a probability space
$(\Omega,\mathcal{A},\P)$. Throughout the paper by a distribution
function we will mean its left-continuous version. Denote
$$
F_n(x)=F^{*n}(x\sqrt{n})= \P\big(X_1+\ldots+X_n<x\sqrt{n}\big),
$$
$$
\Phi(x)= \frac1{\sqrt{2\pi}}\int_{-\infty}^x e^{-t^2/2}\,dt,\quad
x\in\R.
$$
The classical Berry--Esseen theorem states that there exists a
finite positive absolute constant $C_0$ which guarantees the
validity of the inequality
$$
\ud\equiv\sup_x|F_n(x)-\Phi(x)|\le C_0\beta_3/\sqrt{n}\eqno(1)
$$
for all $n\ge1$ and any $F\in\F_3$ (Berry, 1941), (Esseen, 1942).
The problem of establishing the best value of the constant $C_0$ in
inequality (1) has a long history and is very rich in deep and
interesting results. A detailed history of the efforts to lower the
upper estimates of $C_0$ from the original works of A.~Berry (Berry,
1941) and C.-G.~Esseen (Esseen, 1942) to the papers of I.~Shiganov
(Shiganov, 1982), (Shiganov, 1986) can be found in (Korolev and
Shevtsova, 2009). Here we will give an outline of the recent history
of the subject.

In 2006 I.~Shevtsova improved Shiganov's upper estimate
$C_0\le0.7655$ by approximately $0.06$ and obtained the estimate
$C_0\le0.7056$ (Shevtsova, 2006). In 2008 she sharpened this
estimate to $C_0\le0.7005$ (Shevtsova, 2008). In 2009 the mutually
beneficial competition for the best estimate of the constant became
especially keen. On 8 June, 2009 I.~Tyurin submitted his paper
(Tyurin, 2010a) to the <<Theory of Probability and Its
Applications>>. That paper, along with other results, contained the
estimate $C_0\le0.5894$. Two days later the summary of those results
was submitted to <<Doklady Akademii Nauk>> (translated into English
as <<Doklady Mathematics>>) (Tyurin, 2009b). Independently, on 14
September, 2009 V.~Korolev and I.~Shevtsova submitted their paper
(Korolev and Shevtsova, 2009) to the <<Theory of Probability and Its
Applications>>. In that paper the inequality
$$
\ud\le 0.34445\cdot \frac{\beta_3+0.489}{\sqrt{n}}, \ \ \ n\ge1,
$$
was proved which holds for any distribution $F\in\F_3$ yielding the
estimate $C_0\le0.5129$ by virtue of the condition $\beta_3\ge1$.

On 17 November, 2009 the paper (Tyurin, 2010c) was submitted to the
<<Russian Mathematical Surveys>> (its English version (Tyurin,
2009d) appeared on 3 December, 2009 on arXiv:0912.0726). In this
paper together with the other results the estimate $C_0\le0.4785$
was proved.

By a modification of the method used in (Korolev and Shevtsova,
2009) with sharpened estimate of the difference of characteristic
functions in the vicinity of zero the authors of the mentioned work
proved the inequality
$$
\ud\le0.33477\cdot\frac{\beta_3+0.429}{\sqrt{n}}\eqno(2)
$$
in the paper (Korolev and Shevtsova, 2010) submitted to
<<Scandinavian Actuarial Journal>> on 16 March 2010 and published
online on 04 June, 2010. Inequality (2) also improves the upper
bound for the absolute constant $C_0$, since it implies the estimate
$C_0\le0.33477\cdot1.429<0.4784$.

In his oral communication at the 10th International Vilnius
Conference on Probability Theory and Mathematical Statistics on 1
July, 2010 I.\,Tyurin announced the same estimate as in (Korolev and
Shevtsova, 2010): $C_0\le0.4784$. However, in the corresponding
abstract (Tyurin, 2010) no concrete values of the constant $C_0$
were presented.

In his talk at the 3rd Northern Triangular Seminar organized by the
Euler International Mathematical Institute (an international office
of St.-Petersburg Department of Steklov Institute of Mathematics) on
11 April, 2011, the presentation (Tyurin, 2011) of which is
available at
http://www.pdmi.ras.ru/EIMI/2011/NTS/presentations/tyurin.pdf,
I.\,Tyurin announced the estimate $C_0\le0.4774$. However, the proof
of that result was not given in (Tyurin, 2011).

As concerns the lower estimates for $C_0$, in 1956 C.-G. Esseen
found the bound: $C_0\ge C_E$ with
$$
C_{E}=\frac{\sqrt{10}+3}{6\sqrt{2\pi}}=0.409732\ldots
$$
(Esseen, 1956). In 1967 V.~Zolotarev put forward the hypothesis
that $C_0=C_E$ in (1) (Zolotarev, 1967a), (Zolotarev, 1967b).
However, up till now this hypothesis has been neither proved nor
rejected.

By a modification of the method used in (Korolev and Shevtsova,
2010), here we prove the following

\begin{theorem} For all $n\ge1$ and all $F\in\F_3$ we have
$$
\ud\le 0.3328\cdot\frac{\beta_3+0.429}{\sqrt{n}}.\eqno(3)
$$
\end{theorem}

Obviously, inequality (3) sharpens (2) for all $F\in\F_3$ and
$n\ge1$. Moreover, under the conditions imposed on the moments of
the random variable $X_1$ we always have $\beta_3\ge1$. Therefore,
$0.3328(\beta_3+0.429)\le 0.3328\cdot1.429\beta_3 < 0.4756\beta_3$,
which slightly improves the estimate $C_0\le0.4774$ announced in
(Tyurin, 2011). However, the method proposed below allows to improve
the upper bound for $C_0$ still more, if replace $0.429$ by another
constant. Namely, the following statement holds.

\begin{theorem}For all $n\ge1$ and all $F\in\F_3$ we have
$$
\ud\le 0.33554\cdot\frac{\beta_3+0.415}{\sqrt{n}}.\eqno(4)
$$
\end{theorem}

\begin{remark}
Inequality (4) is more precise than (3) for $\beta_3<1.2854\ldots\
.$ For other values of~$\beta_3$ inequality~(3) is better.
\end{remark}

Since $0.33554(\beta_3+0.415)\le 0.33554\cdot1.415\beta_3 <
0.4748\beta_3$, we obtain

\begin{corollary} The classical Berry--Esseen inequality $(1)$ holds
with $C_0=0.4748$.
\end{corollary}

\begin{remark}
Inequality (4) is sharper than the classical Berry--Esseen
inequality (1) with the best known constant $C_0=0.4748$ presented
in corollary 1 for all possible values of $\beta_3\ge1$. Corollary 1
also slightly improves the estimate of the absolute constant in the
classical Berry--Esseen inequality announced in (Tyurin, 2011).
\end{remark}

\begin{remark}
Even if the hypothesis of V.\,M.~Zolotarev that
$C_0=C_{E}=0.4097\ldots$ in (1) (see (Zolotarev, 1967a), (Zolotarev,
1967b)) turns out to be true, then, due to that $\beta_3\ge1$,
inequalities (3) and (4) will be sharper than the classical
Berry--Esseen inequality (1) for $\beta_3\ge 1.86$ and $\beta_3\ge
1.88$ respectively.
\end{remark}

\section{Proofs}

Here we use the approach proposed and developed by V.\,M.~Zolotarev
in his works (Zolotarev, 1965), (Zolotarev, 1966), (Zolotarev,
1967a), (Zolotarev, 1967b), modified in (Korolev and Shevtsova,
2010). Below we will point out only the main ideas which distinguish
this work from the previous ones (for details see (Korolev and
Shevtsova, 2010)).

Denote
$$
f(t)=\E e^{itX_1},\quad
f_n(t)=\left(f\left(\frac{t}{\sqrt{n}}\right)\right)^n,\quad r_n(t)
= |f_n(t)-e^{-t^2/2}|,\quad t\in\R.
$$

\begin{lemma}[Prawitz, 1972]\label{PrawitzSmoothInequality}
For an arbitrary distribution function $F$ and $n\ge1$ for any
$0<t_0\le1$ and $T>0$ we have the inequality
$$ \ud\le 2\int_0^{t_0}|K(t)|r_n(Tt)\,dt+
2\int_{t_0}^{1}|K(t)|\cdot|f_n(Tt)|dt+
$$
$$
+2\int_0^{t_0}\left|K(t)-\frac i{2\pi t}\right|e^{-T^2t^2/2}dt +
\frac1{\pi}\int_{t_0}^\infty e^{-T^2t^2/2}\frac{dt}t,
$$
where
$$ K(t)=\frac12(1-|t|)+\frac i2\left[(1-|t|)\cot\pi
t+\frac{\sign t}\pi\right],\quad -1\le t\le1.
$$
\end{lemma}

Now consider the estimates of the characteristic functions appearing
in lemma~\ref{PrawitzSmoothInequality}.

\begin{lemma} For all $t\in\R$
$$
r_n(t)\le 2e^{-t^2/2}\int_0^{|t|}ue^{u^2/2}
\sin\Big(\frac{u\ell}4\wedge\frac\pi2\Big)
\Big|f\Big(\frac{u}{\sqrt{n}}\Big)\Big|^{n-1}du,
$$
$$
r_n(t)\le 2\int_0^{|t|}ue^{u^2/(2n)}
\sin\Big(\frac{u\ell}4\wedge\frac\pi2\Big)du\cdot
\frac1n\sum_{k=0}^{n-1}
\Big|f\Big(\frac{t}{\sqrt{n}}\Big)\Big|^{n-k-1}
\exp\Big\{-\frac{kt^2}{2n}\Big\}.
$$
\end{lemma}

\proof. The first estimate is proved in (Korolev, Shevtsova, 2010),
and the second one follows from the inequality
$$
|a^n-b^n|\le |a-b|\sum_{k=0}^{n-1}|a|^k|b|^{n-k-1},
$$
which is valid for any complex numbers $a$ and $b$, where we put
$a=e^{-t^2/(2n)}$, $b=f(t/\sqrt{n})$, as well as the first estimate
of the lemma with $n=1$ being applied.

\smallskip

Now we proceed to the estimation of $|f(t)|$. For $\eps>0$ and
$t\in\R$ set
$$
\psi(t,\eps)=
\begin{cases}t^2/2-\varkappa\eps|t|^3,&\text{$|t|\le\theta_0/\eps$},\vspace{1mm}\cr
\displaystyle{(1-\cos\eps t)/\eps^2},&\text{$\theta_0<\eps|t|\le
2\pi$},\vspace{1mm}\cr 0,&\text{$|t|>2\pi/\eps$},
\end{cases}
$$
where $\theta_0=3.99589567\ldots$ is the unique root of the equation
$$
\theta^2 + 2\theta\sin \theta + 6(\cos \theta - 1)=0,\quad
\pi\le\theta\le2\pi,
$$
$$
\varkappa \equiv \sup_{x>0}x^{-3}\left|\,\cos
x-1+x^2/2\right|=\theta_0^{-3}\big(\cos \theta_0-1+\theta_0^2/2\big)
=0.09916191\ldots
$$
It can easily be made sure that the function $\psi(t,\eps)$
monotonically decreases in $\eps>0$ for any fixed $t\in\R$. Denote
the Lyapunov fraction by $\ell=\beta_3/\sqrt{n}$.

\begin{lemma} For any $F\in\F_3$, $n\ge1$ the following
estimates hold:
$$
|f_n(t)|\le  \big[1-2\psi(t,\ell+1/\sqrt{n})/n\big]^{n/2} \le
\exp\{-\psi(t,\ell+1/\sqrt{n})\}\le
$$
$$
\le \exp\big\{-t^2/2+\varkappa(\ell+1/\sqrt{n})|t|^3 \big\},\quad
t\in\R,
$$
$$
|f_n(t)|\le \bigg[\Big(1-\frac{\psi(t,\ell)}n\Big)^2+
\frac{\ell^2t^6}{36n^2}\bigg]^{n/2} \le
\exp\Big\{-\psi(t,\ell)+\frac{\psi^2(t,\ell)}{2n}+
\frac{\ell^2t^6}{72n}\Big\},\quad |t|\le\frac\pi2\sqrt{n}.
$$
\end{lemma}

\proof. For the first series of the estimates see (Prawitz, 1973)
and (Korolev, Shevtsova, 2010). Let us prove the second one.
Evidently,
$$
|f_n(t)|=\bigg[\Big(\E\cos\frac{tX_1}{\sqrt{n}}\Big)^2+
\Big(\E\sin\frac{tX_1}{\sqrt{n}}\Big)^2\bigg]^{n/2}.
$$
Since $\E X_1=0$, for all $t\in\R$ we have
$$
\big(\E\sin tX_1\Big)^2=\big(\E(\sin tX_1-tX_1)\Big)^2\le
\Big(\E|\sin tX_1-tX_1|\Big)^2\le
\Big(\frac{|t|^3}{6}\E|X_1|^3\Big)^2= \frac{\beta_3^2|t|^6}{36}.
$$
In (Sakovi\v{c}, 1965) it is proved that for any random variable $X$
with $\E X^2=1$ the inequality $\E\cos tX\ge0$ holds for all
$|t|\le\pi/2$. On the other hand, from (Prawitz, 1973) it follows
that for any random variable $X$ with $\E X^2=1$ and
$\E|X|^3<\infty$
$$
\E\cos(tX)\le 1-\psi(t,\E|X|^3),\quad t\in\R.
$$
Thus, for all $|t|\le\pi\sqrt{n}/2$ we have
$$
\Big|\E\cos \frac{tX_1}{\sqrt{n}}\Big|=
\E\cos\frac{tX_1}{\sqrt{n}}\le 1-\psi(t/\sqrt{n},\beta_3)=
1-\frac{\psi(t,\ell)}{n},
$$
which together with the estimate for $\big(\E\sin tX_1\big)^2$ leads
to the desired result.

\smallskip

Finally, the process of computational optimization can be properly
organized with the help of the following statements.

\begin{lemma}[Bhattacharya and Ranga Rao, 1976] For any distribution $F$
with zero mean and unit variance we have
$$
\rho(F,\Phi)\le \sup_{x>0}\left\{\Phi(x)-x^2/(1+x^2)\right\}=
0.54093654\ldots
$$
\end{lemma}

Moreover, repeating the algorithms described in (Prawitz, 1975) and
(Gaponova, Shevtsova, 2009) we conclude that inequality (3) holds
true at least for $\ell\le0.0357$
and inequality (4) is true at
least for $\ell\le0.0353$. 

The lemmas presented above give the grounds for restricting the
domain of the values of
$\eps=(\beta_3+0.429)/\sqrt{n}=\ell+0.429/\sqrt{n}$ in theorem~1 and
$\eps=(\beta_3+0.415)/\sqrt{n}=\ell+0.415/\sqrt{n}$ in theorem~2 by
a bounded interval separated from zero. The supremum in $n$ is
estimated by using the monotonically decreasing majorants for the
characteristic functions for $n$ large enough. The fact that all the
estimates used for the characteristic functions monotonically
increase in $\eps$ allows to estimate the supremum in $\eps$ by
computation of the quantity under consideration in a finite number
of points as in the preceding works.

\proof\ of theorem 1. Using the algorithm described above we found
one extremal point: $n=4$, $\eps\approx0.8565$
($\beta_3\approx1.284$, $t_0\approx0.398$, $T\approx5.451$), the
extremal value do not exceeding $0.3328$, which proves theorem~1.

\proof\ of theorem 2. The computations carried out according to the
algorithm described above show that there are two extremal points:
$n=4$, $\eps\approx0.838$ ($\beta_3\approx1.261$, $t_0\approx0.394$,
$T\approx5.513$) and $n=6$, $\eps\approx0.5777$ ($\beta_3\approx1$,
$t_0\approx0.317$, $T\approx7.723$). Both extremal values do not
exceed $0.33554$, whence the assertion of theorem~1 follows.

\smallskip

{\bf Acknowledgements.} The author has the pleasure to express her
sincere gratitude to professor V.\,Yu.\,Korolev for permanent
attention to this work.

\bigskip

\small

\centerline{R E F E R E N C E S}

\begin{enumerate}

\item A. C. Berry. The accuracy of the Gaussian approximation to
the sum of independent variables.~-- {\it Transactions of the
American Mathematical Society}, 1941, Vol. 49, p. 122--136.

\item C.-G. Esseen. On the Liapunoff limit of error in the theory
of probability.~-- {\it Ark. Mat. Astron. Fys.}, 1942, Vol. A28, No.
9, p. 1--19.

\item C.-G. Esseen. A moment inequality with an application to the central limit
theorem.~-- {\it Skand. AktuarTidskr.}, 1956, Vol. 39, p. 160--170.

\item M. O. Gaponova and I. G. Shevtsova. Asymptotic estimates of the
absolute constant in the Berry--Esseen inequality for distributions
with infinite third moments.~-- {\it Informatics and Its
Applications}, 2009, Vol. 3, No. 4, p. 41--56 (in Russian).

\item V. Yu. Korolev and I. G. Shevtsova.
On the upper bound for the absolute constant in the Berry-Esseen
inequality.~-- {\it Theory Probab. Appl.}, 2010, vol. 54, No. 4, p.
638--658.

\item V. Yu. Korolev and I. G. Shevtsova. An improvement of the
Berry--Esseen inequality with applications to Poisson and mixed
Poisson random sums.~-- {\it Scandinavian Actuarial Journal}, 04
June 2010. Online first:
http://www.tandfonline.com/doi/abs/10.1080/03461238.2010.485370

\item H.~Prawitz. Limits for a distribution, if the
characteristic function is given in a finite domain.~-- {\it Scand.
AktuarTidskr.}, 1972, p.~138--154.

\item H.~Prawitz. Ungleichungen f\"{u}r den absoluten Betrag einer
charakteristischen funktion.~-- {\it Skand. AktuarTidskr.}, 1973,
No.~1, p.~11--16.

\item H.~Prawitz. On the remainder in the central
limit theorem.~-- {\it Scand. Actuarial J.}, 1975, No.~3,
p.~145--156.

\item G.\,N.\,Sakovi\v{c}. Pro \v{s}irinu spectra.~--
{\it Dopovidi AN URSR}, 1965, No.\,11, p.\,1427--1430.

\item I. G. Shevtsova. A refinement of the upper estimate
of the absolute constant in the Berry--Esseen inequality.~-- {\it
Theory of Probability and its Applications}, 2006, Vol.~51, No.~3,
p.~622--626.

\item I. G. Shevtsova. On the absolute constant in the Berry--Esseen
inequality.~-- In: {\it The Collection of Papers of Young Scientists
of the Faculty of Computational Mathematics and Cybernetics, Moscow
State University}, Issue 5. Publishing House of the Faculty of
Computational Mathematics and Cybernetics, Moscow State University,
Moscow, 2008, p. 101--110 (in Russian).

\item I. G. Shevtsova. An improvement of convergence rate estimates in
Lyapunov's theorem.~-- {\it Doklady Mathematics}, 2010, Vol. 435,
No. 1, p. 26--26.

\item I. S. Shiganov. On a refinement of the upper constant in the
remainder term of the central limit theorem.~-- In: {\it Stability
Problems for Stochastic Models. Proceedings of the Seminar}.
Publishing House of the Institute for Systems Studies, Moscow, 1982,
p. 109--115 (in Russian).

\item I. S. Shiganov. Refinement of the upper bound of the constant
in the central limit theorem.~-- {\it Journal of Soviet
Mathematics}, 1986, Vol. 35, p. 2545--2550.

\item I. S. Tyurin. On the convergence rate in Lyapunov's
theorem.~-- {\it Theory of Probability and its Applications}, 2011,
Vol.\,55, No.\,2, p.\,253--270.

\item I. S. Tyurin. On the accuracy of the Gaussian
approximation.~--  {\it Doklady Mathematics}, 2009, Vol. 80, No. 3,
p. 840--843.

\item I. S. Tyurin. Refinement of the upper bounds of the
constants in Lyapunov's theorem.~-- {\it Russian Mathematical
Surveys}, 2010, Vol.\,65, No.\,3, p.\,586--588.

\item I. Tyurin. New estimates of the convergence rate in the
Lyapunov theorem.~-- arXiv:0912.0726v1, 3 December, 2009.

\item I. S. Tyurin. Some new results concerning the rate of
convergence in Lyapunov’s theorem.~-- {\it Book of Abstracts of 10th
International Vilnius Conference on Probability Theory and
Mathematical Statistics}, 2010, p. 279.

\item I. S. Tyurin. Some new advances in estimating the rate of
convergence in CLT.~-- {\it 3d Northern Triangular Seminar}, 11
April 2011. Presentation of the talk.\\
http://www.pdmi.ras.ru/EIMI/2011/NTS/presentations/tyurin.pdf


\item P. van Beek. An application of Fourier methods to the problem of sharpening
the Berry--Esseen inequality.~-- {\it Z. Wahrsch. verw. Geb.}, 1972,
Bd. 23, S. 187--196.

\item V. M. Zolotarev. On closeness of the distributions of two sums of independent random
variables.~-- {\it Theory of Probability and its Applications},
1965, Vol. 10, No. 3, p. 519--526.

\item V. M. Zolotarev. An absolute estimate of the remainder term
in the central limit theorem.~-- {\it Theory of Probability and its
Applications}, 1966, Vol. 11, No. 1, p. 108--119.

\item V. M. Zolotarev. A sharpening of the inequality of
Berry--Esseen.~-- {\it Z. Wahrsch. verw. Geb.}, 1967, Bd. 8, S.
332--342.

\item V. M. Zolotarev. Some inequalities in probability theory and
their application in sharpening the Lyapunov theorem.~-- {\it Soviet
Math. Dokl.}, 1967, Vol. 8, p. 1427--1430.

\end{enumerate}

\end{document}